\documentclass[reqno,11pt]{amsart}
\usepackage{amssymb,amsmath}

\newtheorem{theorem}{Theorem}[section]

\newtheorem{proposition}[theorem]{Proposition}
\newtheorem{corollary}[theorem]{Corollary}

\theoremstyle{definition}

\theoremstyle{remark}
\newtheorem{remark}[theorem]{Remark}

\numberwithin{equation}{section}

\newcommand{\re}{\mathbb{R}}
\newcommand{\N}{\mathbb{N}}

\def\lap{\Delta}

\def\qed{\hfill\framebox(6,8)\par\vspace{1em}}

\def\weakto{\rightharpoonup}
\def\ol{\overline}

\def\pd{\partial}
\def\la{\lambda}

\def\H01{H_0^1(\Omega)}
\def\intO{\int_{\Omega}}

\def\intR2{\int_{\re^2}}

\def\1{P_1}
\def\2{P_2}
\def\({\left(}
\def\){\right)}

\begin{document}

\title[Two blow up points]{On the location of two blow up points on an annulus for the mean field equation}


%    author one information
%\author{}
%\address{}
%\curraddr{}
%\email{}
%\thanks{}

\author[M. Grossi]{Massimo Grossi}
\address{Dipartimento di Matematica, Universit\`a di Roma
``La Sapienza", P.le A. Moro 2 - 00185 Roma} 
\email{massimo.grossi@uniroma1.it}

\author[F. Takahashi]{Futoshi Takahashi}
\address{Department of Mathematics, Osaka City University \& OCAMI,
Sumiyoshi-ku, Osaka, 558-8585, Japan}
\email{futoshi@sci.osaka-cu.ac.jp}

%    author two information
%\author{}
%\address{}
%\curraddr{}
%\email{}
%\thanks{}

%    \subjclass is required.
\subjclass[2010]{Primary 35J08, Secondary 35J25, 35J60.}
%35B33, % Critical exponents
%35B40, % Asymptotic behavior of solutions
%35J20, % Variational methods for second-order elliptic equations
%35J08, % Green's function
%35J25, % Boundary value problems for second-order elliptic equations
%35J60. % Nonlinear PDE of elliptic type

\date{\today}

%\dedicatory{}

%    "Communicated by" -- provide editor's name; required.
\commby{J. Krieger}

%    Abstract is required.
\begin{abstract}
We consider the mean field equation on two-dimensional annular domains,
and prove that if $\1$ and $\2$ are two blow up points of a blowing-up solution sequence of the equation, 
then we must have $\1=- \2$.
\end{abstract}

\maketitle


%%%%%%%%%%%%%%%%%%%%%%%%%%%%%%%%%%%%%%%%%%%%%%%%%%%%%%%%%%%%%%%%%%%%%%%%%%%%%%%%
%
%   \S 1. Introduction.
%
%%%%%%%%%%%%%%%%%%%%%%%%%%%%%%%%%%%%%%%%%%%%%%%%%%%%%%%%%%%%%%%%%%%%%%%%%%%%%%%%
%\sezione{Introduction}\label{s1}
\section{Introduction}\label{s1}
\vspace{1em}
In this paper we consider the problem
\begin{equation}
\label{MFE}
    \begin{cases}
        -\lap u = \la \frac{e^u}{\intO e^u dx} & \quad \mbox{in} \; \Omega, \\
        u = 0 & \quad \mbox{on} \; \pd\Omega,
    \end{cases}
\end{equation}
where $\Omega$ is a smooth bounded domain in $\re^2$ and $\la > 0$ is a parameter.
The equation (\ref{MFE}) is known as the mean field equation and is considered to have relations with
various fields of mathematical physics, such as
Onsager's vortex theories, Chern-Simons-Higgs gauge theory, and so on.
The interested readers should refer the books by Tarantello \cite{Tarantello(book)}, Yang \cite{Yang(book)}, and the references therein.
The possible blowing-up or non-compactness for a solution sequence of the problem have attracted many authors for more than two decades,
and many efforts have been devoted to study such a critical phenomena.

Now, thanks to the works by \cite{Nagasaki-Suzuki}, \cite{Brezis-Merle} and \cite{Ma-Wei}, 
we have the following description of the blowing-up solution sequences:
Let $u_n$ be a sequence of solutions to (\ref{MFE}) for $\la = \la_n$ such that
$\| u_n \|_{L^{\infty}(\Omega)}$ is not bounded from above while $\la_n = O(1)$ as $n \to \infty$.
Then there exists a subsequence $\la_n$ and 
a set $\mathcal{S} = \{ a_1, \cdots, a_l \}$
with $a_i\in\Omega$,
such that $\la_n \to 8\pi l, \; l \in \N$, and
\[
	\la_n\frac{e^{u_n}}{\intO e^{u_n} dx} \weakto 8\pi \sum_{i=1}^l \delta_{a_i}
\]
in the sense of measures.
Moreover, each $a_i \in \mathcal{S}$ must satisfy the condition
\begin{equation}
\label{Chara}
%    \frac{1}{2} \nabla R(a_i) - \frac{1}{8\pi} \nabla \log V (a_i) - \sum_{j = 1, j \ne i}^l \nabla_x G(a_i, a_j) = \vec{0}.
	\frac{1}{2} \nabla R(a_i) - \sum_{j = 1, j \ne i}^l \nabla_x G(a_i, a_j) = \vec{0}, \quad (i=1,2,\cdots,l)
\end{equation}
where $G = G(x,y)$ is the Green function with pole $y \in \Omega$ subject to the Dirichlet boundary condition:
\[
	-\lap_x G(x,y) = 2 \pi \delta_y \quad \text{in} \; \Omega, \quad G(x,y) \Big|_{x \in \pd\Omega} = 0,
\]
and $R$ is the Robin function defined as
\[
	R(y) = \lim_{x \to y} \( \log |x-y|^{-1} - G(x,y) \). 
\]
Therefore, the relation (\ref{Chara}) can be considered as a 
characterization of the location of blow up points for $(\ref{MFE})$.

On the other hand, several existence results of $l$-points blowing-up solutions to $(\ref{MFE})$ have been found by several authors,
see \cite{Esposito-Grossi-Pistoia}, \cite{delPino-Kowalczyk-Musso}.
Their results can be summarized as follows:

Let $l \ge 1$ be an integer and 
set 
\[
	\Delta = \{(x_1, \cdots, x_l) \in \Omega^l \; | \; x_i = x_j \; \text{for some} \; i,j \in \{1,\cdots,l\} \},
\]
where $\Omega^l \subset \re^{2l}$ denotes an $l$-time products of $\Omega$.
Define $\mathcal{F}: \Omega^l \setminus \Delta \to \re$ as
\[
	\mathcal{F}(\xi_1,\cdots,\xi_l) = \sum_{i=1}^l R(\xi_i) - \sum_{i \ne j \atop 1 \le i,j \le l} G(\xi_i, \xi_j),
\]
here, we agree that $\mathcal{F}(\xi) = R(\xi)$ for $\xi \in \Omega$ when $l = 1$.
Note that the condition $\nabla_{(\xi_1,\cdots,\xi_l)} \mathcal{F}(a_1, \cdots, a_l) = 0$ is equivalent to 
(\ref{Chara}) for $(a_1, \cdots, a_l) \in \Omega^l$.
By these notations,
let $(a_1, \cdots, a_l) \in \Omega^l \setminus \Delta$ be a ``stable" critical point \cite{Esposito-Grossi-Pistoia},
or a ``nontrivial" critical point \cite{delPino-Kowalczyk-Musso} of $\mathcal{F}$,
that is, $(a_1, \cdots, a_l)$ satisfies (\ref{Chara}) and some additional ``stability" or ``nontriviality" condition is satisfied.
Then there exists a sequence of solutions blowing up exactly at $\mathcal{S} = \{ a_1, \cdots,a_l \}$.
In particular, if the domain is not simply-connected, there always exists  a sequence of blowing-up solution
which blows up at $l$ points on the domain for any $l \in \N$.
Contrary to the above,
we do not have any blowing-up solution sequence with multiple $(l \ge 2)$ blow up points, if the domain is convex.
This nonexistence of multiple blow up points holds true for several nonlinear problems other than (\ref{MFE}), see \cite{Grossi-Takahashi}.
The relationship between the location of blow up points and the geometry of the domain seems to be an interesting subject.

In this note, we turn to the study of the location of blow up points for the mean field equation (\ref{MFE}).
We concentrate to the case when $\Omega$ is an annulus.
In this case, C. C. Chen and C. S. Lin \cite{Chen-Lin(AIHP)} showed the following:

\begin{theorem} (\cite{Chen-Lin(AIHP)} Theorem 1.4.)
\label{Theorem:Chen-Lin}
Let $\{ u_n  \}$ be a solution sequence to (\ref{MFE}) for $\la = \la_n$ with $\la_n \to 16 \pi$ such that 
$u_n$ blows up at two points $P_1$ and $P_2$ on the annulus,
Let $P_{1,n}$ and $P_{2,n}$ be the two local maximum points near $P_1$ and $P_2$ respectively, 
then $P_{1,n}$, $P_{2,n}$ and the origin form a straight line $l_n$ and
$u_n$ is symmetric with respect to the line $l_n$ for $n$ large.
Consequently, $P_1, P_2$ and the origin are located on a same line.
\end{theorem}

The proof of Theorem \ref{Theorem:Chen-Lin} is done by the method of rotating planes,
which is applicable to other kinds of nonlinear elliptic equations, see for example \cite{Lin-Takagi}. An analogous result for problems involving the critical Sobolev exponent was obtained in  \cite{CP}.\\
Theorem  \ref{Theorem:Chen-Lin} leaves open the question of whether the blow up points $P_1$ and  $P_2$ are anti-symmetric, i.e.
\begin{equation}\label{i2}
\1=-\2.
\end{equation}
In this note, by using the characterization of blow up points (\ref{Chara}) and the explicit form of the Green function on an annulus
derived by D. M. Hickey \cite{Hickey(Ann.Math)}, \cite{Hickey(Bull.AMS)}, we show \eqref{i2}.
\begin{theorem}
\label{Theorem:Main}
Let $\{ u_n  \}$ be a sequence of solutions to (\ref{MFE}) for $\la = \la_n$ with $\la_n \to 16 \pi$ such that 
$u_n$ blows up at two points $\1$ and $\2$ on the annulus,
Then we have $\1=-\2$.
\end{theorem}
Next we compute  the value of $|\1|=|\2|$.
\begin{theorem}\label{1}
%There is a unique $\sqrt{R R_1} < r_0 < R^{3/4} R_1^{1/4}$ such that $|a| = |b| = r_0$. 
Define $r_0 = |\1| = |\2|$ where $\1, \2 \in A = \{ a< |x| <a \}$ are two blow up points. Then  $r_0$ is the unique solution of the equation
\begin{align}
\label{condition1}
2 \frac{\log (r/b)}{\log (a/b)} - \frac{1}{2} 
= \sum_{m=1}^{\infty} \frac{1}{b^{2m} - a^{2m}} ( r^{2m} - (ab)^{2m} r^{-2m} ) ((-1)^m + 1)
\end{align}
for $r\in(a,b)$.
\end{theorem}
The explicit form of the Dirichlet Green function on a two dimensional annulus can be seen in several literatures,
see for example, \cite{Courant-Hilbert}, \cite{Akhiezer}, \cite{Aksoy-Celebi}. 
Most of them use the Weierstrass doubly periodic functions.
We find that the Fourier expansion of the Green function is convenient to our analysis.
Since the derivation in \cite{Hickey(Ann.Math)} is easy and seems less known, we prove the formula in Appendix for the sake of completeness.

%%%%%%%%%%%%%%%%%%%%%%%%%%%%%%%%%%%%%%%%%%%%%%%%%%%%%%%%%%%%%%%%%%%%%%%%%%%%%%%%
%
%   \S 2. Proof of Theorem \ref{Theorem:Main}.
%
%%%%%%%%%%%%%%%%%%%%%%%%%%%%%%%%%%%%%%%%%%%%%%%%%%%%%%%%%%%%%%%%%%%%%%%%%%%%%%%%
%\sezione{Proof of Theorem \ref{Theorem:Main}.}\label{s2}
\section{Proof of Theorem \ref{Theorem:Main}.}\label{s2}
\vspace{1em}

Let $A = \{ x \in \re^2 \; | \;a< |x| < b \}$ be a two-dimensional annulus.
Then the Green function on $A$ is explicitly written as follows.

\begin{proposition}(Hickey's formula \cite{Hickey(Ann.Math)}) 
\label{Prop:Hickey}
Let $G_A = G_A(x,y)$ be the Green function on $A$ with pole $y \in A$:
\[
	-\lap_x G_A(x,y) = 2 \pi \delta_y \quad \text{in} \; A, \quad G_A(x,y) \Big|_{x \in \pd A} = 0.
\]
Then we have
\begin{align}
\label{Green_A}
	G_A(x,y) &= -\log |x-y| + A_0(y) + B_0(y) \log |x| \notag \\
	&-\sum_{m=1}^{\infty} \frac{1}{m} ( A_m(y) |x|^m + B_m(y) |x|^{-m} ) \cos m(\theta - \theta_y),
\end{align}
where $x = (x_1, x_2) = (|x| \cos \theta, |x| \sin \theta)$, $y = (|y| \cos \theta_y, |y| \sin \theta_y)$,
and
\begin{align}
\label{AmBm}
	A_0(y) &= \log b \frac{\log (a/|y|)}{\log (a/b)}, \quad B_0(y) = \frac{\log (|y|/b)}{\log (a/b)}, \notag \\
	A_m(y) &= \frac{|y|^m - \( \frac{a^2}{|y|} \)^m}{b^{2m} - a^{2m}}, \quad 
	B_m(y) = \frac{a^{2m} \( \( \frac{b^2}{|y|} \)^m - |y|^{m} \)}{b^{2m} - a^{2m}}.
\end{align}
\end{proposition}

As a corollary, we have
\begin{corollary}
\label{Cor:Robin_A}
The Robin function on the annulus $A = \{a< |x| < b \} \subset \re^2$ is 
\begin{align}
\label{Robin_A}
	R_A(y) &:= \lim_{x \to y} \( -\log |x-y| - G_A(x,y) \) \notag \\ 
	&= -\frac{(\log |y| - \log b)^2}{\log (a/b)} - \log b \notag \\ 
	&+ \sum_{m=1}^{\infty} \frac{1}{m} \frac{1}{b^{2m} - a^{2m}} ( |y|^{2m} - 2a^{2m} + (ab)^{2m} |y|^{-2m} ).
\end{align}
\end{corollary}

Note that $R_A$ is a radial function on $A$, as it was stated in \cite{Chen-Lin(AIHP)} (Lemma 3.3).

Also using the fact
\[
	\nabla_x = \frac{x}{r} \frac{\pd}{\pd r} + \frac{x^{\perp}}{r^2} \frac{\pd}{\pd \theta}
\]
where $r = |x|, x^{\perp} = (-x_2, x_1)$ for $x = (x_1, x_2)$, 
we obtain the formula for the gradients of $G_A$ and $R_A$ as follows: 
\begin{corollary}
\label{Cor:D_Green_AD_Robin_A}
We have
\begin{align}
\label{D_Green_A}
	\nabla_x G_A(x,y) &= -\frac{(x-y)}{|x-y|^2} + B_0(y) \frac{x}{|x|^2} \notag \\
	&-\frac{x}{|x|} \sum_{m=1}^{\infty} ( A_m(y) |x|^{m-1} - B_m(y) |x|^{-m-1} ) \cos m(\theta - \theta_y) \notag \\
	&+ \frac{x^{\perp}}{|x|^2} \sum_{m=1}^{\infty} ( A_m(y) |x|^{m} + B_m(y) |x|^{-m} ) \sin m(\theta - \theta_y),
\end{align}
and
\begin{align}
\label{D_Robin_A}
	\frac{1}{2} \nabla R_A(y) 
	&= - \frac{\log (|y|/b)}{\log (a/b)} \frac{y}{|y|^2} \notag \\ 
	&+ \sum_{m=1}^{\infty} \frac{1}{b^{2m} - a^{2m}} ( |y|^{2m-1} - (ab)^{2m} |y|^{-2m-1} ) \frac{y}{|y|}.
\end{align}
\end{corollary}

\vspace{1em}
Now, we prove Theorem \ref{Theorem:Main} by direct calculations.

\vspace{1em}\noindent
{\bf Proof of Theorem \ref{Theorem:Main}.}
Let $\1, \2 \in A$, $\1\ne\2$ be two blow up points for a blowing-up solution sequence $\{ u_n \}$ to (\ref{MFE}).
Since Theorem \ref{Theorem:Chen-Lin} holds, the only thing we have to prove Theorem \ref{Theorem:Main} is that $|\1| = |\2|$. 
For that purpose, we will exploit the characterization of blow up points (\ref{Chara}). 
In this case, it reads that
\begin{align}
\label{char_blow_up}
	\begin{cases}
	&\frac{1}{2} \nabla R_A(\1) = \nabla_x G_A(\1,\2), \\
	&\frac{1}{2} \nabla R_A(\2) = \nabla_x G_A(\2,\1),
	\end{cases}
\end{align}
which implies
\begin{align}
\label{char_blow_up2}
	\begin{cases}
	&\frac{1}{2} \nabla R_A(\1) \cdot \1 = \nabla_x G_A(\1,\2) \cdot \1, \\
	&\frac{1}{2} \nabla R_A(\2) \cdot \2 = \nabla_x G_A(\2,\1) \cdot \2.
	\end{cases}
\end{align}
By using the formulae (\ref{D_Green_A}), (\ref{D_Robin_A}), we can write the equations (\ref{char_blow_up2}) as
\begin{align}
\label{char_blow_up3_1}
	&-B_0(\1) + \sum_{m=1}^{\infty} \frac{1}{b^{2m} - a^{2m}} ( |\1|^{2m} - (ab)^{2m} |\1|^{-2m} ) \notag \\
	&= -\frac{(\1-\2) \cdot \1}{|\1-\2|^2} + B_0(\2) - \sum_{m=1}^{\infty} ( A_m(\2) |\1|^{m} - B_m(\2) |\1|^{-m} ) \cos m(\theta_{\1} - \theta_{\2}),
\end{align}
and
\begin{align}
\label{char_blow_up3_2}
	&-B_0(\2) + \sum_{m=1}^{\infty} \frac{1}{b^{2m} - a^{2m}} ( |\2|^{2m} - (ab)^{2m} |\2|^{-2m} ) \notag \\
	&= -\frac{(\2-\1) \cdot P_2}{|\2-\1|^2} + B_0(\1) - \sum_{m=1}^{\infty} ( A_m(\1) |\2|^{m} - B_m(\1) |\2|^{-m} ) \cos m(\theta_{\1} - \theta_{\2}),
\end{align}
where $\1 = (|\1| \cos \theta_{\1}, |\1| \sin \theta_{\1})$, $\2 = (|\2| \cos \theta_{\2}, |\2| \sin \theta_{\2})$ in polar coordinates.
Inserting (\ref{AmBm}), we have
%\begin{align*}
%	&A_m(b) |a|^m - B_m(b) |a|^{-m} \\
%&= \frac{1}{R^{2m} - R_1^{2m}} \left\{ |a|^m |b|^m - R_1^{2m} |a|^m |b|^{-m} + R_1^{2m} |a|^{-m} |b|^m - (R_1 R)^{2m} |a|^{-m} |b|^{-m} \right\}, \\
%	&A_m(a) |b|^m -B_m(a) |b|^{-m} \\
%&= \frac{1}{R^{2m} - R_1^{2m}} \left\{ |a|^m |b|^m - R_1^{2m} |a|^{-m} |b|^m + R_1^{2m} |a|^m |b|^{-m} - (R_1 R)^{2m} |a|^{-m} |b|^{-m} \right\}.
%\end{align*}
\begin{align*}
	&A_m(\2) |\1|^m - B_m(\2) |\1|^{-m} = \frac{1}{b^{2m} -a^{2m}} \times \\
	& \left\{ |\1|^m |\2|^m - a^{2m} |\1|^m |\2|^{-m} + a^{2m} |\1|^{-m} |\2|^m - (ab)^{2m} |\1|^{-m} |\2|^{-m} \right\}, \\
	&A_m(\1) |\2|^m -B_m(\1) |\2|^{-m} = \frac{1}{b^{2m} - a^{2m}} \times \\
	& \left\{ |\1|^m |\2|^m - a^{2m} |\1|^{-m} |\2|^m + a^{2m} |\1|^m |\2|^{-m} - (ab)^{2m} |\1|^{-m} |\2|^{-m} \right\}.
\end{align*}
Thus, subtracting (\ref{char_blow_up3_2}) from (\ref{char_blow_up3_1}), we have
\begin{align*}
	&\sum_{m=1}^{\infty} \frac{1}{b^{2m} - a^{2m}} ( |\1|^{2m} - (ab)^{2m} |\1|^{-2m} - |\2|^{2m} + (ab)^{2m} |\2|^{-2m} ) \notag \\
	&= \frac{|\2|^2 - |\1|^2}{|\1 - \2|^2} - \sum_{m=1}^{\infty} \frac{2a^{2m}}{b^{2m} -a^{2m}}
( |\1|^{-m} |\2|^m - |\1|^m |\2|^{-m} ) \cos m (\theta_{\1} - \theta_{\2}).
\end{align*}
From this, we obtain
\begin{align}
\label{subtract}
	&\frac{|\2|^2 - |\1|^2}{|\2 - \1|^2} \notag \\
	&=\sum_{m=1}^{\infty} \frac{|\1|^{2m} - |\2|^{2m}}{b^{2m} - a^{2m}} 
	\left\{ 1 + \frac{(ab)^{2m}}{|\1|^{2m} |\2|^{2m}} - \frac{2a^{2m}}{|\1|^m |\2|^m} \cos m (\theta_{\1} - \theta_{\2}) \right\}. 
\end{align}
Concerning the RHS of (\ref{subtract}), we see
\begin{align*}
	&\left\{ 1 + \frac{(ab)^{2m}}{|\1|^{2m} |\2|^{2m}} - \frac{2a^{2m}}{|\1|^m |\2|^m} \cos m (\theta_{\1} - \theta_{\2}) \right\} \\ 
	&\ge 1 + \frac{(ab)^{2m}}{|\1|^{2m} |\2|^{2m}} - \frac{2a^m b^m}{|\1|^m |\2|^m} 
	= \( 1 - \frac{(ab)^{m}}{|\1|^{m} |\2|^{m}} \)^2 \ge 0,
\end{align*}
since $a < b$.
Thus, if $|\1| > |\2|$, $\mbox{LHS of (\ref{subtract})} < 0$ while $\mbox{RHS of (\ref{subtract})} \ge 0$, which is a contradiction.
The case of $|\1| < |\2|$ leads to the same contradiction.
This implies that $|\1| = |\2|$ must hold, which ends the proof of Theorem \ref{Theorem:Main}.
\qed

Now we compute the value of $|\1| = |\2|$.\\
{\bf Proof of Theorem \ref{1}.}
By inserting $\2 = -\1$ into the first equation of (\ref{char_blow_up}):
\[
	\frac{1}{2} \nabla R_A(\1) = \nabla_x G_A(\1,\2),
\]
and using (\ref{D_Green_A}), (\ref{D_Robin_A}), we have 
\begin{align*}
	&-\frac{\log (|\1|/b)}{\log (a/b)} \frac{\1}{|\1|^2}
+ \frac{\1}{|\1|^2} \sum_{m=1}^{\infty} \frac{1}{b^{2m} - a^{2m}} ( |\1|^{2m} - (ab)^{2m} |\1|^{-2m} ) \notag \\
	&= -\frac{1}{2} \frac{\1}{|\1|^2} +  \frac{\log (|\1|/b)}{\log (a/b)} \frac{\1}{|\1|^2} \\
&- \frac{\1}{|\1|^2} \sum_{m=1}^{\infty} \frac{(-1)^m}{b^{2m} - a^{2m}} ( |\1|^{2m} - (ab)^{2m} |\1|^{-2m} ),
\end{align*}
which in turn implies
\begin{align}
\label{condition1}
	2 \frac{\log (|\1|/b)}{\log (a/b)} - \frac{1}{2} 
= \sum_{m=1}^{\infty} \frac{1}{b^{2m} - a^{2m}} ( |\1|^{2m} - (ab)^{2m} |\1|^{-2m} ) \{ (-1)^m + 1 \}
\end{align}
since $\1 \ne 0$.
Let $f(r) = 2 \frac{\log (r/b)}{\log (a/b)} - \frac{1}{2}$ for $a < r < b$.
$f$ is a monotonically decreasing function with $f(a + 0) = \frac{3}{2}, f(b-0) = -\frac{1}{2}$, and having a unique zero
at $r = b^{3/4} a^{1/4}$.
Also define
\[
	g(r) = \sum_{m=1}^{\infty} \frac{1}{b^{2m} - a^{2m}} ( r^{2m} - (ab)^{2m} r^{-2m} ) \{ (-1)^m + 1 \}.
\]
Since $(-1)^m + 1 \ge 0$ for any $m \in \N$, 
we see $g$ is monotonically increasing with respect to $r$ and
\begin{align*}
&\lim_{r \downarrow a} g(r) = \sum_{m=1}^{\infty} \frac{1}{b^{2m} - a^{2m}} (a^{2m} - b^{2m} ) \{ (-1)^m + 1 \} = -\infty, \\
&\lim_{r \uparrow b} g(r) = \sum_{m=1}^{\infty} \frac{1}{b^{2m} - a^{2m}} ( b^{2m} - a^{2m} ) \{ (-1)^m + 1 \} = +\infty,
\end{align*}
with having unique zero $r = \sqrt{ab}$.
Thus we have the unique $r_0, \sqrt{ab} < r_0 < b^{3/4} a^{1/4}$ such that $f(r_0) = g(r_0)$ by the Intermediate Value Theorem for continuous functions.
\qed
\begin{remark}
By the proof of the last theorem it follows that $\sqrt{ab} < r_0 < b^{3/4} a^{1/4}$.
\end{remark}
\begin{remark}
It is interesting to know  what will happen when the number of blow up points is three or more: 
Up to now, we do not obtain the possible conclusion $|\1| = |\2| = |P_3|$ from the identities
\[
	\begin{cases}
	&\frac{1}{2} \nabla R_A(\1) = \nabla_x G_A(\1,\2) + \nabla_x G_A(\1,P_3), \\
	&\frac{1}{2} \nabla R_A(\2) = \nabla_x G_A(\2,\1) + \nabla_x G_A(\2,P_3), \\
	&\frac{1}{2} \nabla R_A(P_3) = \nabla_x G_A(P_3,\1) + \nabla_x G_A(P_3,\2).
	\end{cases}
\]
We conjecture that if we have $m$-blow up points on the two-dimensional annulus, 
then they must be located on the vertices of regular $m$-polygon. 
The verification of this seems difficult.
\end{remark}

%%%%%%%%%%%%%%%%%%%%%%%%%%%%%%%%%%%%%%%%%%%%%%%%%%%%%%%%%%%%%%%%%%%%%%%%%%%%%%%%
%
%   \S 3. Appendix. Proof of Proposition \ref{Prop:Hickey}
%
%%%%%%%%%%%%%%%%%%%%%%%%%%%%%%%%%%%%%%%%%%%%%%%%%%%%%%%%%%%%%%%%%%%%%%%%%%%%%%%%
%\sezione{Appendix. Proof of Proposition \ref{Prop:Hickey}}
\section{Appendix. Proof of Proposition \ref{Prop:Hickey}}
\label{Appendix}
\vspace{1em}
In this Appendix, we prove Proposition \ref{Prop:Hickey}.
Let $A = \{ x \in \re^2 \; | \;  a < |x| < b \}$ be an annulus in $\re^2$ as before and set
\[
	G_A(x, y) = -\log |x-y| + u(x, y),
\]
where $u(x,y)$ is harmonic with respect to $x \in A$ and coincides with $\log |x-y|$ when $x \in \pd A$.
We use the polar coordinate for $x, y \in \ol{A}$ and write $x = (|x| \cos \theta, |x| \sin \theta)$, $y = (|y| \cos \theta_y, |y| \sin \theta_y)$. 
Take $u(x,y)$ in the form
\begin{align*}
%	&u(x,y) = a_0(y) + b_0(y) \log |x| \\
%	&+ \sum_{m=1}^\infty \left\{ \( a_m(y) |x|^m + b_m(y) |x|^{-m} \) \cos m \theta + \( c_m(y) |x|^m + d_m(y) |x|^{-m} \) \sin m \theta \right\}
	u(x,y) &= a_0(y) + b_0(y) \log |x| \\
	&+ \sum_{m=1}^\infty \( a_m(y) |x|^m + b_m(y) |x|^{-m} \) \cos m \theta \\ 
	&+ \sum_{m=1}^\infty \( c_m(y) |x|^m + d_m(y) |x|^{-m} \) \sin m \theta
\end{align*}
which is harmonic in $x$.
Recalling the expansion
\begin{align*}
	&\log |x-y| \Big|_{|x| = b} = \log b - \sum_{m=1}^{\infty} \frac{1}{m} \( \frac{|y|}{b} \)^m \cos m (\theta - \theta_y) \\
	&= \log b - \sum_{m=1}^{\infty} \frac{1}{m} \( \frac{|y|}{b} \)^m \cos m \theta_y \cos m \theta  
	- \sum_{m=1}^{\infty} \frac{1}{m} \( \frac{|y|}{b} \)^m \sin m \theta_y \sin m \theta,
\end{align*}
for $y \in A$, $x = (b\cos \theta,b\sin \theta)$,
and
\begin{align*}
	&\log |x-y| \Big|_{|x| = a} = \log |y| - \sum_{m=1}^{\infty} \frac{1}{m} \( \frac{a}{|y|} \)^m \cos m (\theta - \theta_y) \\
	&= \log |y| - \sum_{m=1}^{\infty} \frac{1}{m} \( \frac{a}{|y|} \)^m \cos m \theta_y \cos m \theta
	- \sum_{m=1}^{\infty} \frac{1}{m} \( \frac{a}{|y|} \)^m \sin m \theta_y \sin m \theta
\end{align*}
for $y \in A$, $x = (a\cos \theta,a\sin \theta)$.
Thus the conditions
\[
	u(x,y) = \log |x-y| \quad \text{for} \; |x| = b \; \text{and} \; |x| = a
\]
reduces to that
\begin{align*}
	&a_0(y) + b_0(y) \log b = \log b, \\ 
	&a_m(y) b^m + b_m(y) b^{-m} = -\frac{1}{m} \(\frac{|y|}{b}\)^m \cos m \theta_y, \\
	&c_m(y) b^m + d_m(y) b^{-m} = -\frac{1}{m} \(\frac{|y|}{b}\)^m \sin m \theta_y, \\
	&a_0(y) + b_0(y) \log a = \log |y|, \\ 
	&a_m(y) a^m + b_m(y) a^{-m} = -\frac{1}{m} \(\frac{a}{|y|}\)^m \cos m \theta_y, \\
	&c_m(y) a^m + d_m(y) a^{-m} = -\frac{1}{m} \(\frac{a}{|y|}\)^m \sin m \theta_y,
\end{align*}
which in turn implies 
\begin{align*}
	&a_0(y) = A_0(y) = \log b \frac{\log \frac{a}{|y|}}{\log \frac{a}{b}}, \\
	&b_0(y) = B_0(y) = \frac{\log \frac{b}{|y|}}{\log \frac{b}{a}}, \\
	&a_m(y) = -\frac{1}{m} A_m(y) \cos m \theta_y = -\frac{1}{m} \frac{|y|^m - \(\frac{a^2}{|y|}\)^m}{b^{2m} - a^{2m}} \cos m \theta_y, \\
	&b_m(y) = -\frac{1}{m} B_m(y) \cos m \theta_y = -\frac{1}{m} \frac{a^{2m} \( \(\frac{b^2}{|y|}\)^m - |y|^m \)}{b^{2m} - a^{2m}} \cos m \theta_y, \\
	&c_m(y) = -\frac{1}{m} A_m(y) \sin m \theta_y = -\frac{1}{m} \frac{|y|^m - \(\frac{a^2}{|y|}\)^m}{b^{2m} - a^{2m}} \sin m \theta_y, \\
	&d_m(y) = -\frac{1}{m} B_m(y) \sin m \theta_y = -\frac{1}{m} \frac{a^{2m} \( \(\frac{b^2}{|y|}\)^m - |y|^m \)}{b^{2m} - a^{2m}} \sin m \theta_y.
\end{align*}
Thus the Green function on $A = \{ a < |x| < b \}$ is written in the form (\ref{Green_A}).
\qed

\vspace{1em}\noindent
{\bf Acknowledgement.}
Part of this work was done when the second author (F.T) was visiting Universit\`a di Roma in April, 2010.
He thanks Dipartimento di Matematica for its support and  hospitality. 
F.T was also supported by JSPS Grant-in-Aid for Scientific Research (B), No. 23340038,
and JSPS Grant-in-Aid for Challenging Exploratory Research, No. 24654043.

\bibliographystyle{amsplain}

\end{document}